\documentclass[10pt,a4paper]{amsart}

\usepackage{amsmath,amssymb,amsfonts,amsthm,enumerate,color,array}

\newtheorem{lemma}{Lemma}

\newtheorem*{thmcite}{Theorem A}
\theoremstyle{remark}
\theoremstyle{plain}\newtheorem{theorem}{Theorem}

\title[A random graph with duplications and deletions]
{Further properties of a random graph with duplications and
  deletions}
\author{\'Agnes Backhausz}
\address{MTA Alfr\'ed R\'enyi Institute of Mathematics, Budapest and
  Department of Probability Theory and Statistics\\ 
E\"otv\"os Lor\'and University\\
P\'azm\'any P.~s.\ 1/C, H-1117
  Budapest, Hungary}  
 \email{agnes@math.elte.hu}
 \author{Tam\'as F.~M\'ori}
 \address{Department of Probability Theory and Statistics\\
 E\"otv\"os Lor\'and University\\
P\'azm\'any P.~s. 1/C, 
 H-1117 Budapest, Hungary}
 \email{mori@math.elte.hu}
 \dedicatory{\upshape {\sc MTA Alfr\'ed R\'enyi Institute of
  Mathematics, Budapest and  Department of Probability Theory and 
  Statistics}, \\ 
{\sc E\"otv\"os Lor\'and University}\\ 
P\'azm\'any P.~s. 1/C, H-1117 Budapest, Hungary\\                        
 \textit{E-mail address:} \texttt{agnes@math.elte.hu}\\ 
{\sc Department of Probability Theory and
 Statistics}, \\ 
{\sc E\"otv\"os Lor\'and University}\\ 
P\'azm\'any P.~s. 1/C, H-1117 Budapest, Hungary\\                        
\textit{E-mail address:} \texttt{mori@math.elte.hu} }
\thanks{Supported by the Hungarian Scientific Research Fund -- OTKA K 108615.}


\date{18 September 2014}

\begin{document}

\begin{abstract}
We deal with a random graph model where at each step, a vertex is
chosen uniformly at random, and  it is either duplicated or its edges
are deleted. Duplication has a given probability.  We analyse the limit 
distribution of the degree of a fixed vertex, and derive
a.s.\ asymptotic bounds for the maximal degree. The
model shows a phase transition phenomenon with respect to the
probabilities of duplication and deletion.  
\end{abstract}

\maketitle
\thispagestyle{empty}


\section{Introduction}

For researchers in mathematical biology it is evident that duplication
of the information in the genome is a dominant evolutionary force in
shaping biological networks. On the other hand, due to injuries, deletion 
of edges or vertices is also a phenomenon which is natural to
consider. We analyse the random graph model that was described in a
recent paper of Th\"ornblad \cite[2014+]{Th}. This model evolves  
in discrete time steps, and it has a parameter $0<\theta<1$. 
 We start from a single vertex without edges. At each step, we choose a 
vertex $v$ uniformly at random. With probability $\theta$ we
duplicate $v$; that is, we add a new vertex and  
connect it to the neighbours of $v$ and to $v$ itself with single
edges. Otherwise, with probability $1-\theta$, all edges  
of $v$ are deleted (the vertex itself stays in the graph, and has the
chance to get new edges later on).  

As it was presented in \cite{Th} and as we will see later, the
asymptotic behaviour of the model depends on the value of $\theta$,
and there is  a phase transition phenomenon. Naturally,  our
results will also be different for certain regimes of the duplication
probability. The case $0<\theta<1/2$ is the subcritical case, where
deletion is  
more likely, and as we will prove, the maximal degree has the order of
logarithm of the actual number of vertices, almost surely. In the
critical case ($\theta=1/2$) the maximal degree grows faster; we
have the square of  
the logarithm of the number of vertices. Finally, in the supercritical
case, when   $1/2<\theta<1$ and the duplication is dominant, the
maximal degree will be compared to the number of vertices (without
logarithm) in some sense. We remark that similar phase transition is
present in some other random graph models where deletion (but no
duplication) is introduced; see e.g.\ the work of Vallier
\cite[2013]{vallier}.  

The random graph defined above consists of disjoint cliques. The
vertices form separate clusters, and  two vertices are  
connected if and only if they belong to the same cluster. Due to this
simple structure, it may be interesting from  
the point of view of coagulation--fragmentation models or population
dynamics. Champagnat, Lambert and Richard described certain properties
of a continuous version of this model, the so-called  
splitting trees with mutations, where the clusters correspond to
different alleles of a gene \cite[2012]{amaury1, amaury3},
\cite[2013]{amaury2}. 

This shows that the current model does not have  a fine structure as a
graph. In the paper \cite[2014+]{dup} we described the connection of
the critical case of  
the model analysed here to a model which is also built on duplication
and deletion, but has a nontrivial graph structure. That model also
turns out to be highly clustered:  
as one can expect,  due to the duplication steps, some 
dense clusters evolve, while edges between clusters are rare. These
highly clustered networks come up in mathematical biology, e.g.\ for
modelling protein--protein interaction networks. Thus graph models
with duplication (but with a richer graph structure than the actual
one) may also be found in the literature.  
In most of those models edges are never deleted, but only some
randomly chosen edges of the chosen vertex are duplicated, and some
extra random edges may be added to the new vertex (Kim at al.
\cite[2002]{kim}, Pastor-Satorras, Smith and Sol\'e
\cite[2003]{pastor}, Chung et al. \cite[2003]{chung}, Bebek et
al. \cite[2006]{bebek}). However, these 
papers did not contain mathematically rigorous arguments, and some
results of the earlier ones (stating that the degree distribution is
polynomially decaying with an exponential cutoff) were disclaimed by
the latter ones. Recently, Hermann and Pfaffelhuber
\cite[2014+]{pfaff} have proved several results on the frequency of
isolated vertices and cliques, and also on the evolution of the degree
of a fixed vertex in the initial graph. Various other models were
also introduced, where the 
choice of the duplicated vertex is not uniform but depends on the
degrees  (Jordan \cite[2011]{jordan},  Cohen,  Jordan and
Voliotis \cite[2010]{cohen}, Farczadi and Wormald
\cite[2014+]{farczadi}) or on the state of a hidden Markov chain (Hamdi,
Krishnamurthy and Yin \cite[2013+]{hamdi}).  

Notice  that in this  duplication--deletion model vertices of larger
degree are more likely 
to increase their degree, because the probability that one of their
neighbours is duplicated is larger. On the other hand, the probabilty
of decreasing their degree is also larger. However, this is a kind of
a preferential attachment phenomenon. Preferential attachment models
are still popular for modelling web graphs or biological networks 
since the seminal paper of Albert and Barab\'asi \cite[1999]{ba}, and
also from a theoretical point of view. It is  worth mentioning the
model free approaches of Ostroumova, Ryabchenko and Samosvat
\cite[2013]{ostr} and 
Dereich and Ortgiese \cite[2014]{marcel}. However, for example, due to
the possibility that the degree of a vertex can decrease to $1$ in a
single step by deletion, our model does not fit into those frameworks,
in which the degree of a fixed vertex can not decrease.
 Maximal degree was also investigated in certain
preferential attachment models, see e.g. \cite[M\'ori,
2005]{M05}. Bubeck, Mossel and R\'acz showed that the seed graph may
have influence on the limiting distribution of the maximal degree in
some kind of preferential attachment models. \cite[2014]{bubeck}.
Since all vertices of  the initial configuration are deleted after
finitely many steps, this phenomenon does not occur in the current model. 

Preferential attachment models are often investigated due to their
scale free property: the proportion of vertices of degree $d$ tends to
some contant $c_d$ in some sense, and $c_d$ decays polynomially as
$d\rightarrow \infty$. The current model has a rather different
asymptotic degree distribution, as the results of Th\"ornblad
\cite{Th} and \cite{dup} shows. Since we will  use it later on, we sum
up these results, as follows. In accordance with \cite{Th} we introduce
\[
\beta=\frac{\theta}{2\theta-1}\,,\quad \gamma=\frac{1-\theta}{\theta}\,.
\] 

\begin{thmcite}\cite[B--M, 2014+, for $\theta=1/2$]{dup}
  \cite[Th\"ornblad, 2014+, for $\theta\neq 1/2$]{Th}.
Let $X[n,k]$ denote the proportion of vertices of degree $k$ after $n$
steps in   the duplication--deletion model defined above. Then we have 
\[
X[n,k]\rightarrow c_k \quad \text{almost surely as } n\rightarrow
\infty,
\]
where $(c_k)_{k=1}^{\infty}$ is the unique sequence of positive numbers
satisfying the following equations.  
\[
c_0=\frac{1-\theta}{1+\theta}(1+c_1); \qquad
c_k=\frac{k+1}{k+1+\theta}\big(\theta c_{k-1}+(1-\theta)c_{k+1}\big)
\quad (k\geq 1).
\]  
Furthermore, $\sum_{k=0}^{\infty} c_k=1$.
\begin{itemize}
\item Subcritical case. If $0<\theta<1/2$, then 
\[
c_k=\gamma^{-k-1}\int_0^1\frac{t^{k+1}(1-t)^{-1-\beta}}
{(1-\gamma^{-1}t)^{1-\beta}}\,dt\quad (k\geq 0), \text{ and }
\]
\[
c_k\sim(-\beta)^{-1}(1-\beta)^{-\beta}\,\Gamma(1-\beta)\gamma^{-k}
k^{\beta}\quad\text{ as } k\rightarrow\infty.
\]
\item Critical case. For  $\theta=1/2$ we have 
\[
c_k=(k+1)\int_0^\infty\frac{t^k e^{-t}}{(1+t)^{k+2}}\,dt \quad (k\geq
0), \text{ and }
\]
\[
c_k\sim (e\pi)^{1/2}\,k^{1/4}\,e^{-2\sqrt{k}}\quad \text{ as }
k\rightarrow\infty.
\]
\item Supercritical case. If $1/2<\theta<1$, then 
\[
c_k=\gamma\int_0^1\frac{t^{k+1}(1-t)^{\beta-1}}{(1-\gamma t)^{1+\beta}}\,dt
\quad (k\geq 0), \text{ and }
\]
\[
c_k\sim\gamma\,\beta^{\beta}\,\Gamma(\beta+1)k^{-\beta}\quad \text{ as }
k\rightarrow\infty.
\]
\end{itemize}
\end{thmcite}
(Note that our $c_k$ corresponds to Th\"ornblad's $d_{k+1}$.)

Therefore the asymptotic degree distribution
decays exponentially in the subcritical case, polynomially in the
supercritical case, and slower than exponential but faster than
polynomial in the critical case. 

The paper is built up as follows. In Section \ref{sth} we describe
useful variants of the model and analyse the evolution of the number
of vertices.  Section \ref{fixed} contains our results on the asymptotic
behaviour of the degree of a fixed vertex. This is used in Section
\ref{maximal}, where we give bounds for the maximal degree, which are
valid with probability $1$. It will follow that the index of the vertex
with maximal degree tends to infinity, that is, there is no persistent hub
in this model.

\section{Variants of Th\"ornblad's model}
\label{sth}

This is a discrete time model. Let us start from a single vertex. The
graph is modified in two ways: at every step a vertex is selected at
random, with equal probability, then this vertex is either duplicated
or deleted. Duplication means that a new vertex is added to the graph,
and it is connected to the selected vertex and its neighbours. 
Deletion means that the edges of the selected vertex get deleted, but
not the vertex itself. Every step, independently of the
past, is a duplication with probability $\theta$, and a deletion with
probability $1-\theta$ $(0<\theta<1)$. This model is called Version 1
or Th\"ornblad's model \cite{Th}.

Further versions differ from Version 1 only by time transforms. 

In Version 2 the development of the graph is slowed down. Let
$N_{n-1}$ be the number of vertices of the graph after the $(n-1)$th
step. At the $n$th step, the graph is not modified 
with probability $1-N_{n-1}/n$. Otherwise each existing vertex has
equal probability to be selected; the selected vertex (exactly one) is
duplicated or deleted with probabilities $\theta$ and $1-\theta$,
resp.  All the above randomizations are independent of each other. As
we will see, in this way the graph does not change in the majority 
of the steps. 

Version 3 is defined in continuous time. Every vertex is given two
clocks at its birth, which alarm according to independent homogeneous
Poisson processes with rates $\theta$ and $1-\theta$. When the first
clocks rings, the vertex gets duplicated, while the second clock
determines deletion times. In this model, steps occur at an
exponentially accelerating pace.  

Focusing only on moments when something happens all versions look
identical. This makes it possible to choose the most convenient
version for different proofs. 

In all versions, at every moment, the graph is a disjoint union of
complete gaphs (cliques). This is obviously true in the beginning, and
it is easy to see that neither duplication, nor deletion can break this
property.

\begin{lemma}\label{size}
Let us denote the number of vertices after $n$ steps by $N_n$, and in
Version 3, the size of the graph at time $t$ will be denoted by
$N(t)$. Then a.s.
\begin{align*}
N_n&\sim \theta n &&\text{in Version $1$};&&\\
N_n&\sim \zeta\,n^{\theta}&& \text{in Version $2$, where $\zeta$ is a
positive random variable};&&\\ 
N(t)&\sim \eta\,e^{\theta t}&& \text{in Version $3$, where $\eta$ is a
positive random variable}.&&
\end{align*}
\end{lemma}
\textbf{Proof.} The proof for Version 1 is obvious.

Version 2. Let  $\mathcal F_n$  denote the $\sigma$-field generated by the
first $n$ steps. At step $n$ the number of vertices increases by $1$ with
(conditional) probability $\theta N_{n-1}/n$. Thus 
\[
E\bigl(N_n\bigm\vert\mathcal F_{n-1}\bigr)=\Big(1+\frac{\theta}{n}
\Big)N_{n-1}.
\]
Introduce 
\[
\varkappa_n=\prod_{i=1}^n \Big(1+\frac{\theta}{i}\Big)^{-1}=
\frac{\Gamma(1+\theta)\Gamma(n+1)}{\Gamma(n+1+\theta)}
\sim \Gamma(1+\theta)n^{-\theta}.
\]
Then $(\varkappa_n N_n,\,\mathcal F_n)$ is a nonnegative martingale,
which is known to be a.s. convergent. Let
$\zeta=\lim_{n\to\infty}n^{-\theta}N_n$. 

We still have to show that $\zeta>0$. Let $R_n=(N_n-1)^{-1}$, if
$N_n\ge 2$, and $R_n=1$ otherwise. This time let
\[
\kappa_n=\prod_{i=1}^n\left(1-\frac{\theta}{i}I(N_{i-1}\ge
2)\right)^{\!-1}\asymp n^\theta 
\]
as $n\to\infty$. Here $I(\,\cdot\,)$ stands for the indicator of the event
in brackets. Clearly,
\[
E\bigl(R_n\bigm|\mathcal F_{n-1}\bigr)=\frac{1}{N_{n-1}}\cdot\frac{\theta
N_{n-1}}{n}+\frac{1}{N_{n-1}-1}\Big(1-\frac{\theta N_{n-1}}{n}\Big)
=\frac{1}{N_{n-1}-1}\Big(1-\frac{\theta}{n}\Big)
\]
on the event $\{N_{n-1}\ge 2\}\in\mathcal F_{n-1}$, and $=1$ on its
complement. Hence $(\kappa_nR_n,\,\mathcal F_n)$ is a nonnegative
martingale. Consequently, $n^{\theta}/N_n$ converges a.s., and its
limit is obviously $1/\zeta$.

Version 3. $N(\theta^{-1})$ is a Yule process (see e.g. \cite{karlin, yule}),
thus it is geometrically distributed, namely $\mathrm{Geom}
\big(e^{-\theta t}\big)$, and $N(t)\sim\eta e^{\theta t}$,  where
$\eta$ is an exponential random variable of expected value $1$. \qed

\section{The degree process of a fixed vertex}
\label{fixed}

In this section we consider Version 3, because, as we will see, that
is the most natural choice for the individual degree
processes. Of course, the results of this section can easily be
transferred to Version 1, by using Lemma \ref{size}.

Let the vertices be labelled by $1,\,2,\,\dots$ in the order of
birth, and let $d_i(t)$ denote the degree of vertex $i$ at time $t$. 

\begin{theorem}\label{degstac}
For every $i=1,2,\dots$ we have
\[
\lim_{t\to\infty}P\big(d_i(t)=k\big)=q_k,\quad k=0,\,1,\,\dots,
\]
where 
\begin{equation}\label{stacmego}
q_0=\gamma(1-c_0),\quad q_k=c_{k-1}-\gamma c_k,\ k=1,2,\dots,
\end{equation}
where the sequence $(c_k)$ is defined in Theorem A. Here
$q_k>0$, $k=0,1,\dots$, and $\sum_{k=0}^\infty q_k=1$.
\end{theorem}
\textbf{Proof.}
Let us fix a vertex. 
Clearly, its degree is a continuous time Markov process with
infinitesimal generator  
\[
\mathbf{M}=\left\lceil
\begin{matrix}
1-\theta&0&\dots\\
1-\theta&0&\dots\\
\vdots&\vdots
\end{matrix}\right\rceil
+\left\lceil
\begin{matrix}
-1 &  \theta \\
1-\theta & -2 &  2\theta\\
 &  2(1-\theta) & -3 &  3\theta\\
  &   &  3(1-\theta) & -4 &  4\theta \\
 &  & &\ddots&\ddots&\ddots
\end{matrix}\right\rceil
\]
The process is positive recurrent, because deletion cuts back the
degree to $0$ at a constant rate. Hence it has a stationary distribution 
$q=(q_0,\,q_1,\,\dots)$ which is the unique discrete distribution
satisfying $q\mathbf M=0$ \cite{karlin}.  Thus,
\begin{equation}\label{stac}
q_0=(1-\theta)(1+q_1),\quad
q_k=\frac{k\theta q_{k-1}+(k+1)(1-\theta)q_{k+1}}{k+1}\,,\ k\ge 1.
\end{equation}

From Theorem A it follows that the numbers $q_k$ in
\eqref{stacmego} satisfy \eqref{stac}. They sum up to $1$, because
\[
\sum_{k=0}^\infty q_k=\gamma(1-c_0)+\sum_{k=1}^\infty(c_{k-1}-
\gamma c_k)=\gamma+(1-\gamma)\sum_{k=0}^\infty c_k=1.
\]
Finally, their positivity follows from the integral form of $c_k$,
which can be found in \cite{dup} for the critical case, and in
\cite{Th} for the subcritical and supercritical cases; see Theorem
A. Namely, we immediately obtain in the subcritical case
\[
q_k=\gamma^{-k}\int_0^1\frac{t^k
(1-t)^{-\beta}}{(1-\gamma^{-1}t)^{1-\beta}}\,dt;
\]
and in the supercritical case
\[
q_k=\gamma\int_0^1\frac{t^k(1-t)^{\beta-1}}
{(1-\gamma t)^\beta}\,dt
\]
for $k\ge 1$. In the critical case, by equation \eqref{stac} and
partial integration we get  
\begin{align*} 
q_k&=k\int_0^\infty\frac{t^{k-1} e^{-t}}{(1+t)^{k+1}}\,dt-
(k+1)\int_0^\infty\frac{t^k e^{-t}}{(1+t)^{k+2}}\,dt\\
&=k\int_0^\infty\frac{t^{k-1} e^{-t}}{(1+t)^{k+1}}\,dt+
\left[\frac{t^k e^{-t}}{(1+t)^{k+1}}\right]_0^\infty
-\int_0^\infty\frac{kt^{k-1} e^{-t}-t^k e^{-t}}{(1+t)^{k+1}}\,dt\\
&=\int_0^\infty\frac{t^k e^{-t}}{(1+t)^{k+1}}\,dt.\qed
\end{align*}

It is somewhat surprising. In spite that the limit distribution of
the degree is the same \textit{for every vertex}, the asymptotic
degree distribution of the graph is different. If the the degrees were
independent and identically distributed, the proportion of vertices
with fixed degree would converge to the corresponding probability. In
our model neither condition is satisfied, not even approximately. On
one hand, in Version 3, the size of the graph grows exponentially. 
Consequently, at every moment the vast majority of the vertices are 
relatively young, so the limit distribution cannot be applied to
them. On the other hand, if the vertices were
nearly independent, the number of vertices with high degree would
follow the Poisson distribution; but in our model, if there exists at
least one vertex of a large degree $d$, then all its neighbours have
the same degree, therefore many of such vertices are coexistent. 

This phenomenon can be better understood if we consider the degree
process of an arbitrary vertex. The higher the degree is, the shorter
it sustains. Therefore a reversed size biased sampling can be
observed: at a given moment the probability that a given vertex has
degree $d$ is less than the proportion of degree $d$ ones among
all vertices.

In Section \ref{maximal} we shall need asymptotics for the tail
of the stationary distribution. 

\begin{theorem}\label{stacasymp}
\begin{align}
\label{substac}
&\text{Subcritical case.}&&q_k+q_{k+1}+\dots\sim
(1-\beta)^{2-\beta}\Gamma(1-\beta)\,k^{\beta-1}\gamma^{-k}\,;&&\\
\label{critstac}
&\text{Critical case.}&&q_k+q_{k+1}+\dots\sim
(e\pi)^{1/2}\,k^{1/4}\,e^{-2\sqrt{k}}\,;&&\\
\label{supstac}
&\text{Supercritical case.}&&q_k+q_{k+1}+\dots\sim
\gamma\,\beta^{\beta}\Gamma(\beta-1)\,k^{1-\beta},&&
\end{align}
as $k\to\infty$.
\end{theorem} 
\textbf{Proof.}
In the subcritical case we have
\[
\sum_{j=k}^\infty q_j=\gamma^{-k}\int_0^1\frac{t^k(1-t)^{-\beta}}
{(1-\gamma^{-1}t)^{2-\beta}}\,dt,
\]
while in the supercritical case
\[
\sum_{j=k}^\infty q_j=\gamma\int_0^1\frac{t^k(1-t)^{-2+\beta}}
{(1-\gamma t)^\beta}\,dt.
\]
If $k\to\infty$, both integrals become relatively negligible over any interval
$[0,\,1-\varepsilon]$, compared to those over
$[1-\varepsilon,\,1]$. Hence in the denominators we can replace $t$ 
with $1$, thus reducing to complete beta integrals. 

In the critical case $q_k+q_{k+1}+\dots=c_{k-1}$, hence Theorem A
can be applied.\qed

Obviously, every vertex becomes isolated infinitely many times, due to
deletion. What can be said about the extremely high degrees?

\begin{theorem}\label{vertmax}
\begin{align*}
&\text{Subcritical case.}&&
\limsup_{t\to\infty}\frac{d_i(t)}{\log\log N(t)}=\frac{1}{\log\gamma}
\,;&&\\
&\text{Critical case.}&&
\limsup_{t\to\infty}\frac{d_i(t)}{(\log\log N(t))^2}=1,&&\\
&\text{Supercritical case.}&&
\limsup_{t\to\infty}\frac{\log d_i(t)}{\log\log N(t)}=
\frac{1}{\beta-1}\,.&&
\end{align*}
for $i=1,2,\dots$ .
\end{theorem}
\textbf{Proof.}
First we investigate how large can the degree grow
between two consecutive deletions. Let $p_i(r)$ denote the probability
that a vertex of degree $i$  will sometimes have degree $r$ at least
once before it is selected for deletion. Then we clearly have
\begin{align*}
p_0(r)&=\theta p_1(r);\\
p_i(r)&=\frac{(i+1)\theta p_{i+1}(r)+i(1-\theta)p_{i-1}(r)}{i+1}\,,\quad
i=1,2,\dots,r-1;\\
p_r(r)&=1.
\end{align*}
Introduce $a_i=p_i(r)/p_0(r)$; it does not depend
on $r$ provided $r>i$. The probability we are interested in is
$p_0(r)=1/a_r$. The sequence $(a_r)_{r\ge 0}$ satisfies the
following recursion.
\begin{equation}\label{arec}
a_0=1,\quad a_1=\frac{1}{\theta}\,,\quad a_i=\theta
a_{i+1}+\frac{i}{i+1}\,(1-\theta)a_{i-1},\quad 
i=1,2,\dots\ .
\end{equation}
In the critical case we can use some well-known facts about Laguerre
polynomials $L_r(x)$ \cite{Szego}. They can be defined by the following
recursion formula.
\begin{multline}\label{legrec}
L_0(x)=1,\ L_1(x)=1-x,\\ 
L_{r+1}(x)=\frac{(2r+1-x)L_r(x)-rL_{r-1}(x)}{r+1}\,,\ r=1,2,\dots\ .
\end{multline}
Their asymptotic behaviour for large $r$ and fixed $y>0$ is given by
\begin{equation}\label{laguerre}
L_r(-y)\sim
2^{-1}\pi^{-1/2}r^{-1/4}e^{-y/2}y^{-1/4}e^{2\sqrt{yr}}.
\end{equation}
Recursions \eqref{arec} and \eqref{legrec} coincide if $\theta=1/2$
and $x=-1$. Hence we obtain 
\begin{equation}\label{critepoch}
p_0(r)=\frac{1}{a_r}=\frac{1}{L_r(-1)}\sim 2\sqrt{e\pi}\,r^{1/4}e^{-2\sqrt{r}}
\end{equation}
in the critical case.

If $\theta\ne 1/2$, we can analyse the asymptotic behaviour of the
sequence $(a_r)$ by computing its generating function
$G(z)=\sum_{r=0}^\infty a_rz^r$. From \eqref{arec} it follows that
\[
\sum_{r=1}^\infty (r+1)a_rz^r=(1-\theta)\sum_{r=1}^\infty ra_{r-1}z^r
+\theta\sum_{r=1}^\infty(r+1)a_{r+1}z^r,
\]
that is,
\[
\big(zG(z)\big)'-1=(1-\theta)z\big(zG(z)\big)'+
\theta\Big(G'(z)-\frac{1}{\theta}\Big).
\]
This leads to the following homogeneous linear ODE.
\[
\big(\theta-z+(1-\theta)z^2\big)G'(z)=\big(1-(1-\theta)z\big)G(z),
\quad G(0)=1.
\]
Its solution can easily be expanded into a power series.
\[
G(z)=(1-z)^{-\beta}(1-\gamma z)^{\beta-1}=
\sum_{r=0}^\infty (-z)^r\sum_{i=0}^r\binom{\beta-1}{i}
\binom{-\beta}{r-i}\gamma^i.
\]
Thus,
\begin{equation*}
a_r=(-1)^r\sum_{i=0}^r\binom{\beta-1}{i}
\binom{-\beta}{r-i}\gamma^i.
\end{equation*}

Suppose first that $\theta>1/2$. Then $\gamma<1$ and $\beta>1$. Since
\[
(-1)^r r^{1-\beta}\binom{-\beta}{r}=\frac{r^{1-\beta}\Gamma(r+\beta)}
{\Gamma(r+1)\Gamma(\beta)}
\] 
converges as $r\to\infty$, therefore it is bounded. Consequently, we
have  
\[
\left|(-1)^r r^{1-\beta}\binom{\beta-1}{i}\binom{-\beta}{r-i}
\gamma^i\right|\le b_i,
\]
uniformly in $r\ge i$, where the infinite series $\sum b_i$ converges. Hence 
\begin{multline*}
\lim_{r\to\infty}r^{1-\beta}a_r=
\sum_{i=0}^\infty \lim_{r\to\infty}(-1)^r r^{1-\beta}\binom{\beta-1}{i}
\binom{-\beta}{r-i}\gamma^i\\
=\frac{1}{\Gamma(\beta)}
\sum_{i=0}^\infty \binom{\beta-1}{i}(-\gamma)^i=\frac{
(1-\gamma)^{\beta-1}}{\Gamma(\beta)}=\frac{\beta^{1-\beta}}{\Gamma(\beta)}\,,
\end{multline*}
and, by this,
\begin{equation}\label{supepoch}
p_0(r)\sim\beta^{\beta-1}\Gamma(\beta)r^{1-\beta}
\end{equation}
in the supercritical case.

The subcritical case is easy to reduce to the supercritical
one. Let $a_r'=\gamma^{-r} a_r$. Then $a_r'$ satisfies the same recursion
that $a_r$ does when $\theta$ is replaced by $1-\theta$. This
substitution transforms the subcritical case into the supercritical
one, furthermore, $\beta$ changes to $1-\beta$. Hence we get
\begin{equation}\label{subepoch}
p_0(r)=\gamma^{-r}(a_r')^{-1}\sim(1-\beta)^{-\beta}\Gamma(1-\beta)\,
r^{\beta}\gamma^{-r}.
\end{equation}   

Up to time $t$ there are $(1-\theta)t(1+o(1))\sim \gamma\log N(t)$
epochs (time intervals between 
consecutive deletions), hence $\max\{d_i(s):s\le t\}$ is
asymptotically equal to the maximum of $(1+o(1))\gamma\log N(t)$
i.i.d. random variables with distribution $P(\xi\ge
r)=p_0(r)$. Starting from \eqref{subepoch}, \eqref{critepoch} and
\eqref{supepoch}, standard 
Borel--Cantelli arguments yield
\[
\max\{d_i(s):s\le t\}\sim\frac{\log\log N(t)}{\log\gamma}
\]
in the subcritical case,
\[
\max\{d_i(s):s\le t\}\sim\log^2t\sim(\log\log N(t))^2
\]
in the critical case, and
\[
\log\max\{d_i(s):s\le t\}\sim\frac{\log\log N(t)}{\beta-1}
\]
in the supercritical case, completing the proof. (Alternatively, one
can apply \cite[Theorem 4.4.4]{Galambos}.)\qed

\section{Maximal degree}

\label{maximal}

Let $M_n$ denote the maximal degree in Version 1 after $n$ steps. 
From Theorem A it is clear that$M_n\to\infty$.  In many scale-free random graph
  processes the order of magnitude of the maximal degree $M_n$ can be
  characterized in the following way: $M_n\asymp \min\{d: N_n c_d<1\}$,
  where $(c_d)$ is the asymptotic degree distribution, and $N_n$ is
  the size of the graph, see e.g.  \cite[2005]{M05}, \cite[2007]{M07}, \cite[2010]{M09},
  \cite[2014]{harom}. This would give $M_n\asymp\log N_n$ in the subcritical case,
$M_n\asymp\log^2 N_n$ in the critical one, and $M_n\asymp
N_n^{1/\beta}$ in the supercritical one. We will show that this
estimate is valid in the subcritical and critical cases, but in the
supercritical case we can prove less. 

\begin{theorem}\label{maxup}
\begin{align*}
&\text{Subcritical case.}&&\frac{1-\theta}{\log\gamma}\le
\liminf_{n\to\infty}\frac{M_n}{\log N_n}\le
\limsup_{n\to\infty}\frac{M_n}{\log N_n}
\le\frac{1+\theta}{\theta\log\gamma}\,.\\
&\text{Critical case.}&&\frac{1}{16}\le\liminf_{n\to\infty}
\frac{M_n}{\log^2N_n}\le\limsup_{n\to\infty}\frac{M_n}{\log^2 N_n}
\le\frac{9}{4}\,.\\
&\text{Supercritical case.}&& \frac{\theta}{\beta}\le
\liminf_{n\to\infty}\frac{\log M_n}{\log N_n}\le\limsup_{n\to\infty}
\frac{\log M_n}{\log N_n}\le\frac{1}{\beta}\,.
\end{align*}

\end{theorem}

\textbf{Proof of the upper bounds.}

The proof will be given for Version 2. 

Let $d_i(n)$ denote the degree of vertex $i$ after step $n$,
$i=1,\dots, N_n$, where $N_n$ is the size of the graph after $n$
steps. Introduce
\[
S_n(r)=\sum_{i=1}^{N_n}\binom{d_i(n)}{r}.
\]
\begin{lemma}\label{esn2}
For every $n=1,2,\dots$ and $r=0,1,2,\dots $ we have
\begin{align}
\label{subbound} &\text{Subcritical case.}&&
ES_n(r)\le 2(r+1)(-\beta)^r n^{\theta};&&\\
\label{critbound} &\text{Critical case.}&&
ES_n(r)\le 2(r+1)!\sqrt{n}\,;&&\\
\label{supbound} &\text{Supercritical case.}&&
ES_n(r)\le\left\{
\begin{array}{ll}
C_r(\theta)\,n^\theta,&\text{if } r<\beta-1;
\rule[-6pt]{0pt}{18pt}\\
C_r(\theta)\,n^\theta(1+\log n),&\text{if }r=\beta-1;
\rule[-6pt]{0pt}{18pt}\\
C_r(\theta)\,n^{(r+1)(2\theta-1)},&\text{if }r>\beta-1.
\rule[-6pt]{0pt}{18pt}\\
\end{array}
\right.
\end{align}
where $C_r(\theta)$ is a constant depending on $r$ and $\theta$ but
independent of $n$..
\end{lemma}
\textbf{Proof.}
First we will verify the following recursion. For every $n=1,2,\dots$ and
$r=1,2,\dots$ we have 
\begin{equation}\label{esnrec}
ES_n(r)=\Big(1+\frac{(2\theta-1)(r+1)}{n}\Big)ES_{n-1}(r)+
\frac{\theta(r+1)}{n}\,ES_{n-1}(r-1).
\end{equation}

At the $n$th step the $i$th term of $S_n(r)$ can
change in the following way. With the notation $d=d_i(n-1)$, 
\[
\binom{d_i(n)}{r}=
\left\{
\begin{array}{cll}
\dbinom{d+1}{r}&\text{with conditional probability}&
\dfrac{d+1}{n}\,\theta\\
&\multicolumn{2}{l}{\text{(vertex $i$ or one of its neighbours is
    duplicated);}}\\ 
\dbinom{d-1}{r}&\text{with conditional probability}&
\dfrac{d}{n}(1-\theta)\\
&\multicolumn{2}{l}{\text{(a neighbour of vertex $i$ is deleted);}}\\
0&\text{with conditional probability}& \dfrac{1-\theta}{n}\\
&\multicolumn{2}{l}{\text{(vertex $i$ is deleted);}}\\ 
\dbinom{d}{r}&\text{otherwise.}&
\end{array}
\right.
\] 
Thus, 
\[
E\biggl(\binom{d_i(n)}{r}\biggm|\mathcal F_{n-1}\biggr)
=\binom{d}{r}\left(1-\frac{d+1}{n}\right)
+\binom{d+1}{r}\frac{d+1}{n}\,\theta
+\binom{d-1}{r}\frac{d}{n}(1-\theta).
\]
Besides, when vertex $i$ is duplicated, an additional term
$\binom{d+1}{r}$ also appears as the yield of the new vertex.
Hence the total contribution of vertex $i$ in
$E\bigl(S_n(r)\bigm|\mathcal F_{n-1}\bigr)$ is
\begin{multline*}
\binom{d}{r}\left(1-\frac{d+1}{n}\right)
+\binom{d+1}{r}\frac{d+2}{n}\,\theta
+\binom{d-1}{r}\frac{d}{n}(1-\theta)\\
=\binom{d}{r}\left(1+\frac{(2\theta-1)(r+1)}{n}\right)
+\binom{d}{r-1}\frac{\theta(r+1)}{n}\,.
\end{multline*}
This implies that
\[
E\bigl(S_n(r)\bigm|\mathcal F_{n-1}\bigr)=
\left(1+\frac{(2\theta-1)(r+1)}{n}\right)S_{n-1}(r)+
\frac{\theta(r+1)}{n}\,S_{n-1}(r-1),
\]
as needed.

Next, we prove the lemma by double induction over $r$ and $n$, basing
on the recursion \eqref{esnrec}.

Clearly, $S_n(0)=N_n$. From the proof of Lemma \ref{size} we know that
\[
EN_n=\prod_{i=1}^n\Big(1+\frac{\theta}{i}\Big)\le 2
\prod_{i=1}^{n-1}\Big(1+\frac{\theta}{i}\Big)\le
2\prod_{i=1}^{n-1}\Big(1+\frac 1i\Big)^{\theta}=2n^\theta
\]
in all three cases. Furthermore, $ES_1(1)=2\theta$ and $ES_1(r)=0$ if
$r>1$. Thus, for all pairs $(n,0)$ and $(1,r)$ Lemma \ref{esn2} holds true.

Let us check the induction step.

In the subcritical case, by using the induction hypothesis we can write 
\begin{multline*}
ES_n(r)\leq 2(r+1)(-\beta)^r(n-1)^\theta\Big(1+\frac{(r+1)\theta}{n\beta}
\Big)\\
\qquad+2r(-\beta)^{r-1}(n-1)^\theta\,\frac{\theta(r+1)}{n}
\le 2(r+1)(-\beta)^{r}\,n^\theta,
\end{multline*}
as needed.

In the critical case we have
\begin{multline*}
ES_n(r)=ES_{n-1}(r)+\frac{r+1}{2n}\,ES_{n-1}(r-1)\\
\le 2(r+1)!\,\sqrt{n-1}\,\Big(1+\frac{1}{2n}\Big)\le
2(r+1)!\,\sqrt{n}.
\end{multline*}

Finally, in the supercritical case $C_0(\theta)=2$ will do. Suppose we
have proved inequality \eqref{supbound} for $r-1$ (and all
$n$). Introduce $s=(r+1)(2\theta-1)$ and
\[
\kappa_n=\frac{\Gamma(n+1+s)}{\Gamma(n+1)}\,,
\]
then
\[
\frac{\kappa_n}{\kappa_{n-1}}=1+\frac sn\le\Big(\frac
{n}{n-1}\Big)^s,
\]
because $1+\frac sn <\big(1+\frac 1n\big)^s\le
\big(1+\frac{1}{n-1}\big)^s=\big(\frac{n}{n-1}\big)^s$, if $s\ge 1$,
and $1+\frac sn\le\big(1-\frac sn\big)^{-1}\le\big(1-\frac 1n\Big)^{-s}
=\big(\frac{n}{n-1}\big)^s$, if $0<s<1$. Hence, for $1\le j\le n$ we
have  
\begin{equation}\label{kappa}
\frac{\kappa_n}{\kappa_j}\le\Big(\frac nj\Big)^s.
\end{equation}
 
By iterating equation \eqref{esnrec} we get
\begin{multline*}
\frac{ES_n(r)}{\kappa_n}=\frac{ES_{n-1}(r)}{\kappa_{n-1}}+
\frac{(r+1)\theta}{n\kappa_n}\,ES_{n-1}(r-1)=\dots\\
=\frac{ES_1(r)}{\kappa_1}+(r+1)\theta\sum_{j=2}^n
\frac{ES_{j-1}(r-1)}{j\kappa_j}\,,
\end{multline*}
hence, by \eqref{kappa},
\[
ES_n(r)\le n^s ES_1(r)
+(r+1)\theta\,n^s\sum_{j=1}^{n-1}j^{-s}ES_j(r-1)=A+B.
\]
In the right-hand side $A$ vanishes if $r>1$. For $r=1$ it
is equal to $2\theta n^{2(2\theta-1)}$, which satisfies
\eqref{supbound} in all three cases. Let us turn to $B$.

First, suppose $r<\beta-1$, that is, $\theta>s$. Then,
by the induction hypothesis we have 
\[
B\le (r+1)\theta\,C_{r-1}(\theta)\,n^s\sum_{j=1}^{n-1}
j^{\theta-s-1}
\le \frac{(r+1)\theta}{\theta-s}\,C_{r-1}(\theta)\,n^{\theta}.
\]

Next, let $r=\beta-1$, that is, $\theta=s$. Then again
\[
B\le (r+1)\theta\,C_{r-1}(\theta)\,n^{\theta}\sum_{j=1}^{n-1}\frac 1j
\le
(r+1)\theta\,C_{r-1}(\theta)\,n^{\theta}(1+\log n).
\]

If $\beta-1<r\le\beta$, that is, $r(2\theta-1)\le\theta<s$, then
\[
B\le (r+1)\theta\,C_{r-1}(\theta)\,n^s\sum_{j=1}^{n-1}
j^{\theta-s-1}(1+\log j)
\le (r+1)\theta\,C_{r-1}(\theta)\,Q\,n^s,
\]
where
\[
Q=\sum_{j=1}^\infty j^{\theta-s-1}(1+\log j)<\infty.
\]

Finally, if $\beta<r$, then
\[
B\le (r+1)\theta\,C_{r-1}(\theta)\,n^s\sum_{j=1}^{n-1}
j^{-2\theta}\le
(r+1)\theta\,C_{r-1}(\theta)\,\zeta(2\theta)\,n^s,
\]
where $\zeta(\,.\,)$ is the Riemann zeta function.\qed

Let us continue the proof of Theorem \ref{maxup}.

In the subcritical case, let us fix $z$ and $a$ in such a way that
$0<z<-1/\beta$, and $a\log(1+z)>1+\theta$. Then, by Lemma \ref{esn2} we have 
\begin{multline*}
E\bigg(\sum_{i=1}^{N_n}(1+z)^{d_i(n)}\bigg)
=E\bigg(\sum_{i=1}^{N_n}\sum_{r=0}^n\binom{d_i(n)}{r}z^r\bigg)
=E\bigg(\sum_{r=0}^n\sum_{i=1}^{N_n}\binom{d_i(n)}{r}z^r\bigg)\\
=\sum_{r=0}^n ES_n(r)z^r\le 2\sum_{r=0}^\infty(r+1)(-\beta z)^rn^\theta
=\frac{2n^\theta}{(1+\beta z)^2}=K\,n^\theta.
\end{multline*}

By the Markov inequality,
\begin{multline*}
P(M_n\ge a\log n)= P\left((1+z)^{M_n}\ge (1+z)^{a\log n}\right)\\
\le n^{-a\log(1+z)}\,E\bigg(\sum_{i=1}^{N_n}(1+z)^{d_i(n)}\bigg)\le
K\,n^{\theta-a\log(1+z)}.
\end{multline*}
The infinite sum of the right-hand side  is convergent as $n$ runs through the positive
integers, thus the Borel--Cantelli lemma implies 
$M_n<a\log n$ a.s. for every sufficiently large $n$. Consequently,\
\[
\limsup_{n\to\infty}\frac{M_n}{\log n}\le 
\frac{1+\theta}{\log(1-\beta^{-1})}=\frac{1+\theta}{\log\gamma}\,.
\]
From Lemma \ref{size} we know that $\log N_n\sim\theta\log n$ as
$n\to\infty$. 

In the critical case we can make use of Laguerre polynomials
again. Their explicit form is
\begin{equation*}
L_k(y)=\sum_{r=0}^k\binom kr\frac{(-y)^r}{r!}\,.
\end{equation*}

Since the
multiplicity of the maximal degree is at least $M_n+1$, we have  
\[
S_n(r-1)\ge(M_n+1)\binom{M_n}{r-1}=r\binom{M_n+1}{r}.
\]
Therefore, by Lemma \ref{esn2},
\[
E\binom{M_n}{r}\le\frac 1r\,ES_n(r-1)\le 2r!\sqrt{n}\,,
\]
for $r=1,2,\dots\ $, hence 
\[
E\big(L_{M_n}(-y)\big)=E\bigg(\sum_{r=0}^n\binom{M_n}{r}\frac{y^r}{r!}
\bigg)\le 1+2\sqrt{n}\sum_{r=1}^n y^r\le\frac{2\sqrt{n}}{1-y}
\]
if $0<y<1$. Let $k=k_n\ge a\log^2 n$, where $ya>9/16$. Then the Markov
inequality, combined with \eqref{laguerre}, implies  
\begin{multline*}
P(M_n\ge k)=P\big(L_{M_n}(-y)\ge L_k(-y)\big)\le \frac
{E\big(L_{M_n}(-y)\big)}{L_k(-y)}\\
=O\Big(\sqrt{n}\,e^{-2(1+o(1))\sqrt
{yk}}\Big)=O\big(n^{1/2-2\sqrt{ya}}\big).
\end{multline*}
The exponent in $O(\,.\,)$ is less than $-1$, hence it makes a
convergent series again, and from the Borel--Cantelli lemma
\[
\limsup_{n\to\infty}\frac{M_n}{\log^2n}\le a
\]
follows for every $a>9/16$. This time $\log n\sim 2\log N_n$ by Lemma
\ref{size}. 

Finally, let us turn to the supercritical case. Let $a>2\theta-1$ and $r$ so
large that $r>\beta-1$ and $r(2\theta-1-a)<-1$ hold. Then by the Markov
inequality and \eqref{supbound} we have
\begin{multline*}
P\big(M_n\ge n^a\big)=P\bigg(\binom{M_n}{r}\ge\binom{n^a}{r}\bigg)\le
\frac{E\binom{M_n}{r}}{\binom{n^a}{r}}\\
\le\frac{ES_n(r-1)}{r\binom{n^a}{r}}=O\big(n^{r(2\theta-1-a)}\big).
\end{multline*}

The proof can be completed with the help of the Borel--Cantelli lemma
and Lemma \ref{size}.\qed

\textbf{Proof of the lower bounds.}

The proof will be performed for Version 3. Let $\varepsilon$ be an
arbitrarily small fixed positive number. The proof will consist of the
following steps. 

We first show that at time $(1-\theta)\sqrt{n}$ there are quite many
isolated points in the graph. Clearly, they behave independently of
each other after time $(1-\theta)\sqrt{n}$. 

Then we give a lower bound for the probability that such an
isolated vertex increases its degree above
$k_n$ by  time $\sqrt{n}$, where $k_n$ is an increasing positive sequence
depending on $\theta$. 
It will   follow that the probability that none of them can do it is
so small that its sum over $n$ is convergent. Hence the
Borel--Cantelli lemma implies that a.s. $M\big(\sqrt{n}\big)\ge k_n$
if $n$ is large enough.  

Finally, we will show that the probability that a vertex having such a
high degree at time $t_n$ will lose from its degree at least
$\varepsilon k_n$ times in the interval $[\sqrt{n},\,\sqrt{n+1}]$ is very
small: it is also finitely summable. Thus, with $n=\lfloor
t^2\rfloor$ we have a.s. $M(t)\ge (1-\varepsilon)k_n$ if $t$ is
large enough.

In more details, let us start with the number of isolated
vertices. For the sake of brevity denote
$N\big((1-\theta)\sqrt{n}-1\big)$ by $N_n$. For each vertex 
present at time $(1-\theta)\sqrt{n}-1$, the probability that, during the
next time unit, it will be deleted some time and not duplicated after
that is 
\[
\sigma=(1-\theta)\big(1-e^{-1}\big).
\]
 (Something happens to
the vertex, and the last event is a deletion.) Therefore the number of
isolated vertices at time $(1-\theta)\sqrt{n}$ is at least as big as a
binomial random variable $Q_n$ with 
parameters $N_n$ and $\sigma$. Since $N_n$ itself is $\mathrm{Geom}(p_n)$
distributed with parameter
$p_n=\exp(-\theta(1-\theta)\sqrt{n}+\theta)$, straightforward
calculation gives that the distribution of $Q_n$ is a mixture: 
\[
Q_n=\left\{
\begin{array}{cll}
\mathrm{Geom}\bigg(\dfrac{p_n}{\sigma+(1-\sigma)p_n}\bigg)
&\text{ with weight }&\dfrac{\sigma}{\sigma+(1-\sigma)p_n}\,;\\
0&\text{ with weight }&\dfrac{(1-\sigma)p_n}{\sigma+(1-\sigma)p_n}\,.
\end{array}
\right.
\]

Let us turn to the estimation of the probability that a fixed isolated
vertex can considerably increase its degree in a time interval of length
$\theta\sqrt{n}$. 

How fast is the convergence to the stationary distribution? This can
be answered easily by coupling. Let us start two degree processes, one
from the stationary distribution, and another one from an isolated
vertex (i.e., from degree $=0$). Let the deletion of the vertex in
question be governed by the same Poisson process in both cases. After
the first deletion stick the two processes together. Then the
probability that the two processes differ after time $t$ is at most
$e^{-(1-\theta)t}$. Hence the same bound is valid for the total variation
distance of the degree distribution at time $t$ from the stationary
one.

Consequently, if a vertex is isolated at time $(1-\theta)\sqrt{n}$,
then the probability that its degree at time $\sqrt{n}$ is larger than
$k_n$, is at least
\[
\pi_n=\sum_{k=k_n}^\infty q_k-\exp\big(-\theta(1-\theta)\sqrt{n}\big).
\]
If $k_n$ is specified in such a way that
\begin{equation}\label{suppose1}
\sum_{n=1}^\infty\frac{p_n}{\pi_n}<\infty
\end{equation}
holds, then $p_n=o(\pi_n)$, and we have
\begin{multline*}
P\left(M\big(\sqrt{n}\big)\le k_n\right)\le E\left((1-\pi_n)^{Q_n}\right)\\
=\frac{(1-\sigma)p_n}{\sigma+(1-\sigma)p_n}+\frac{\sigma}
{\sigma+(1-\sigma)p_n}\cdot\frac{\dfrac{p_n}{\sigma+(1-\sigma)p_n}
(1-\pi_n)}{1-\dfrac{\sigma(1-p_n)}{\sigma+(1-\sigma)p_n}(1-\pi_n)}
\sim\frac{p_n}{\sigma\pi_n}\,.
\end{multline*}
The sum of these probabilities is convergent by supposition. Hence the
  Borel--Cantelli lemma implies that, almost surely,  $M\big(\sqrt{n}\big)>
  k_n$ if $n$ is large enough.

Finally, we have to show that $M(t)$ cannot decrease significantly between
$\sqrt{n}$ and $\sqrt{n+1}$. Suppose $M\big(\sqrt{n}\big)>k_n$. Choose
a vertex with maximal degree and select $k_n$ from its
neighbours. Let us compute the probability that more than $\varepsilon
k_n$ of them will be deleted between $\sqrt{n}$ and
$\sqrt{n+1}$. The number $Z$ of deleted vertices is 
binomial with parameters $k_n$ and 
\[
1-\exp\big((1-\theta)\big(\sqrt{n}-\sqrt{n+1}\big)\big)\le
(1-\theta)\big(\sqrt{n}-\sqrt{n+1}\big)\le 2(1-\theta)n^{-1/2}.
\]
By Hoeffding's inequality 
\[
P\big(Z>\varepsilon k_n\big)\le \exp\left(-2\big(\varepsilon-
2(1-\theta)n^{-1/2}\big)^2 k_n\right).
\]
If, in addition to \eqref{suppose1}, sequence $(k_n)$
  satisfies 
\begin{equation}\label{suppose2}
\sum_{n=1}^\infty \exp\left(-\varepsilon^2 k_n\right)<\infty,
\end{equation}
then the Borel--Cantelli lemma gives
\[
\min\big\{M(t): \sqrt{n}\le t\le \sqrt{n+1}\big\}>
(1-\varepsilon)k_n
\]
if $n$ is sufficiently large. Consequently, with $n=\lfloor
t^2\rfloor$ we have a.s. $M(t)\ge (1-\varepsilon)k_n$ for all
sufficiently large $t$. 

Let us specify $k_n$ in all three cases to meet conditions
\eqref{suppose1} and \eqref{suppose2}.

In the subcritical case let
\[
k_n=\frac{\theta(1-\theta)}{\log\gamma}(1-\varepsilon)\sqrt{n}.
\]
Then condition \eqref{suppose2} is satisfied. Moreover, by
\eqref{substac} we have 
\[
\pi_n=\exp\big(-(1+o(1))\log\gamma\cdot k_n\big).
\]
Hence 
\[
\frac{p_n}{\pi_n}=\exp\big(-(1+o(1))\,\theta(1-\theta)
\varepsilon \sqrt{n}\big),
\]
thus condition \eqref{suppose1} is satisfied as well. Consequently, 
\[
M(t)>(1+o(1))(1-\varepsilon)^2\frac{\theta(1-\theta)}{\log\gamma}t
=(1+o(1))(1-\varepsilon)^2\frac{1-\theta}{\log\gamma}\log N(t)
\]
if $t$ is sufficiently large.

In the critical case let
\[
k_n=\frac{(1-\varepsilon)^2}{64}\,n,
\]
then \eqref{suppose2} is fulfilled. In addition, \eqref{critstac}
implies 
\[
\pi_n=\exp\left(-(1+o(1))\frac{1-\varepsilon}{4}\sqrt{n}\right).
\]
Therefore
\[
\frac{p_n}{\pi_n}=\exp\left(-(1+o(1))\frac{\varepsilon}{4}\sqrt{n}\right), 
\]
and requirement \eqref{suppose1} is also met. Hence
\[
M(t)>(1+o(1))\frac{(1-\varepsilon)^3}{64}\,t^2=
(1+o(1))\frac{(1-\varepsilon)^3}{16}\,\log^2N(t)
\]
if $t$ is large enough.

Finally, in the supercritical case set
\[
k_n=\exp\big((1-\varepsilon)\theta(2\theta-1)\sqrt{n}\big).
\]
Then \eqref{suppose2} is satisfied. By \eqref{supstac} we can write
\[
\pi_n=k_n^{(1+o(1))(1-\beta)}=\exp\big(-(1+o(1))(1-\varepsilon)\theta
(1-\theta)\sqrt{n}\big),
\]
from which it follows that
\[
\frac{p_n}{\pi_n}=\exp\big(-(1+o(1))\varepsilon\theta(1-\theta)\sqrt{n}
\big).
\]
This produces a finite sum, thus \eqref{suppose1} holds
true. Consequently, with $n=\lfloor t^2\rfloor$,
\begin{multline*}
M(t)>(1-\varepsilon)k_n=\exp\big((1+o(1))(1-\varepsilon)
\theta(2\theta-1)t\big)\\
=\exp\big((1+o(1))(1-\varepsilon)(2\theta-1)\log N(t)\big),
\end{multline*}
if $t$ is sufficiently large.\qed

Due to deletions, in our graph there is no persistent hub in the sense
of Krapivsky and Redner \cite[2001]{krap} or Galashin \cite[2014+]{galashin} (namely, 
a single vertex which emerges forever as vertex of maximal degree),
unlike in certain preferential attachment models \cite[Dereich and M\"orters, 2009]{DM}, \cite[M\'ori, 2005]{M05}. As a
corollary to Theorems \ref{vertmax} and \ref{maxup}, it follows that
the index of the vertex with the maximal degree tends to infinity with time.


\begin{thebibliography}{99}
%

 

%
\bibitem{dup} Backhausz, \'A.,  M\'ori, T.\ F., Asymptotic properties
  of a random graph with duplication. To appear in \textit{J.\ Appl.\
  Probab.} \textbf{52(2)} (2015) 

\bibitem{harom} Backhausz, \'A., M\'ori, T. F., Weights and degrees in
  a random graph model based on $3$-interactions.  
  \textit{Acta Math.\ Hungar.} {\bf 143(1)} (2014), 23--43.

\bibitem {ba} Barab{\'a}si, A-L., and Albert, R., Emergence
  of scaling in random networks, \textit{Science} \textbf{286} (1999),
  509--512.  
%
%
\bibitem{bebek} Bebek, G., Berenbrink, P., Cooper, C., Friedetzky, T.,
  Nadeau, J. and Sahinalp, S.\ C., The degree distribution of the
  generalized duplication model. \textit{Theor.\ Comput.\ Sci.},
  \textbf{369} (2006), 234--249.
%
%

\bibitem{bubeck} Bubeck, S., Mossel, E., R\'acz, M. Z., On the
  influence of the seed graph in the preferential attachment
  model. Preprint. \textit{arXiv:1401.4849}. 

\bibitem{amaury1}N. Champagnat\ and\ A. Lambert, Splitting trees with neutral Poissonian mutations I: Small families, Stochastic Process. Appl. {\bf 122} (2012), no.~3, 1003--1033.

\bibitem{amaury3}N. Champagnat, A. Lambert\ and\ M. Richard, Birth and death processes with neutral mutations, Int. J. Stoch. Anal. {\bf 2012}, Art. ID 569081, 20 pp.

\bibitem{amaury2}N. Champagnat\ and\ A. Lambert, Splitting trees with neutral Poissonian mutations II: Largest and oldest families, Stochastic Process. Appl. {\bf 123} (2013), no.~4, 1368--1414.

\bibitem{chung} Chung, F., Lu, L., Dewey, T.\ G., and Galas, D.\ J.,
  Duplication models for biological networks, \textit{J.\ Comput.\
  Biol.}, \textbf{16} (2003), 677--687.

\bibitem{cohen} Cohen, N., Jordan, J., and Voliotis, M., Preferential
  duplication graphs, \textit{J.\ Appl.\ Probab.}, \textbf{47} (2010),
  572--585.  

%
%
\bibitem{DM} Dereich, S., M\"orters, P., Random networks with
  sublinear preferential attachment: Degree evolutions, \textit{
  Elect.\ J.\ Probab.}, \textbf{14} (2009), Paper no. 43, 1222--1267.

\bibitem{marcel} Dereich, S., Ortgiese, M., Robust analysis of
  preferential attachment models with fitness, \textit{
  Combin.\ Probab.\ Comput.} \textbf{23(3)} (2014), 386--411. 
%
%
%
%
%


%

\bibitem{farczadi}Farczadi, L., Wormald, N., On the degree distribution
  of a growing network model. Preprint. \textit{arXiv:1401.0933}. 

\bibitem{galashin} Galashin, P., Existence of a persistent hub in the
  convex preferential attachment model. Preprint. \textit{arXiv:1310.7513}. 


\bibitem{Galambos} Galambos, J., \textit{The Asymptotic Theory of
    Extreme Order Statistics}, Wiley, New York, 1978.

\bibitem{hamdi} Hamdi, M., Krishnamurthy, V., Yin, G.\ G., Tracking
  the empirical distribution of a Markov-modulated
  duplication-deletion random graph. \textit{arXiv:1303.0050[cs.IT]}. 



\bibitem{pfaff} Hermann, F., Pfaffelhuber, P.,  Large-scale behavior
  of the partial duplication random graph. Preprint. 
  \textit{arXiv:1408.0904}. 

\bibitem{jordan} Jordan, J., Randomised reproducing graphs. 
  \textit{Electron.\ J.\ Probab.}, \textbf{16} (2011), 1549--1562.  

\bibitem{kim} Kim, J., Krapivsky, P.\ L., Kahng, B. and Redner, S.,
  Infinite-order percolation and giant fluctuations in a protein
  interaction network. \textit{Phys.\ Rev.}, E66: 055101(R), 2002. 


\bibitem{krap} Krapivsky, L., Redner, S., Organization of growing random networks.
{\it Phys. Rev. E}
 {\bf 63} (2001), no. 6, 066123.



\bibitem{M05} M{\'o}ri, T.\ F., The maximum degree of the
  Barab\'asi--Albert random tree. \textit{Combin.\ Probab.\ Comput.},
  \textbf{14} (2005), 339--348. 

\bibitem{M07} M{\'o}ri, T.\ F., On a $2$-parameter class of
   scale free random graphs. \textit{Acta Math.\ Hungar.}, \textbf{114}
   (2007), 37--48.

\bibitem{M09} M{\'o}ri, T.\ F., Random multitrees. \textit{Studia
   Sci.\ Math.\ Hungar.}, \textbf{47} (2010), 59--80.

%
%


\bibitem{ostr} Ostroumova, L., Ryabchenko, A., Samosvat, E., 
  Generalized preferential attachment: tunable power-law degree
  distribution and clustering coefficient. In: {\it Algorithms and Models
  for the Web Graph}, Lecture Notes in Computer Science \textbf{8305}
  (2013), 185--202.

\bibitem{pastor} Pastor-Satorras, R., Smith, E., and Sol\'e, R.\ V.,
  Evolving protein interaction networks through gene duplication. 
  \textit{J.\ Theor.\ Biol.}, \textbf{222} (2003), 199--210.


%
%
%
%
%
%
\bibitem{karlin}Karlin, S.,  Pinsky, M. A., {\it An Introduction to
    Stochastic Modeling}, fourth edition, Academic Press, 2010. 
    

\bibitem{Szego}Szeg\H o, G., \textit{Orthogonal Polynomials}, 4th
  edition, Amer.\ Math.\ Soc.\ Colloq.\ Publ., vol. 23, AMS,
  Providence, RI, 1975. 

\bibitem{Th} Th\"ornblad, E., Asymptotic degree distribution of a
  duplication--deletion random graph model. Preprint.
  \textit{arXiv:1408.4268v1}.

\bibitem{yule} Yule, G.\ U., A mathematical theory of evolution, based
  on the conclusions of Dr.\ J.\ C.\ Willis, F.R.S., 
  \textit{Philos.\ Trans.\ R.\ Soc.\ Lond.\ Ser.\ B.}, \textbf{213} 
  (1925), 402--410. 

\bibitem{vallier} T. Vallier, Transition of the degree sequence in the random graph model of Cooper, Frieze, and Vera, Stoch. Models {\bf 29} (2013), no.~3, 341--352. 

%
%
%
%
\end{thebibliography}
\end{document}